\documentclass[12pt]{article} 
 
\usepackage{amsmath,enumerate,amsfonts,amssymb,color,graphicx} 

\begin{document}
\newtheorem{Def}{Definition}[section]
\newtheorem{thm}{Theorem}[section]
\newtheorem{lem}{Lemma}[section]
\newtheorem{rem}{Remark}[section]
\newtheorem{question}{Question}[section]
\newtheorem{prop}{Proposition}[section]
\newtheorem{cor}{Corollary}[section]
\newtheorem{clm}{Claim}[section]
\newtheorem{step}{Step}[section]
\newtheorem{sbsn}{Subsection}[section]
\newtheorem{conj}{Conjection}[section]
\title{A response to 
an article of Xu-Jia Wang} 
\author{YanYan Li \footnote{Department of Mathematics, Rutgers University}~ and Luc Nguyen \footnote{Department of Mathematics, Princeton University}}  
\date{}
\maketitle

\maketitle
\begin{abstract}  This is a response to the article 
arXiv:1212.3130v1 by Xu-Jia Wang, where he attempted to
address a mathematical question we raised.  We point out that,
and explain why,  the article is far from answering our objections.
Moreover, we have more recently found
more serious
trouble in the paper under discussion   based on the  false
assertion there
  that minimal radial functions of superharmonic
functions (with respect to a conformal Laplacian of
a Riemannian metric)
are superharmonic. 
\end{abstract}

\setcounter{section}{0}

\section{}

This is a response to  the article
[W]  by Xu-Jia Wang.    We point out that
the article is far from answering our objections.

Firstly, we  disagree with the claim of Wang,
made in the last paragraph of Section 3 of [W], 
that Section 2 of [W] answers our objections to the proof of the
main theorem  Theorem 1.3 in [TW].
The reasons
 are given
  in the next section.

Secondly, 
 we  don't  understand why
the
equicontinuity of $h^{(r)}(\cdot):=
h(r\cdot)$, for $0<r<1$, when $h$ is locally bounded,
 is  clear to Wang,  as claimed in Section 4 of [W]
(line 4 on page 4).
 Here $h$ is the function
in Lemma 3.4 of [W].
We have repeatedly asked Wang,
starting
from Nov. 16,  to provide a proof of this assertion
 made in his email on Nov. 14, but he has 
 never given one. Pages of detailed
arguments in our paper [LN] can be used to prove this.

Thirdly,
 having studied [TW] in more detail, we have found more serious 
trouble  based on the  false
assertion on line -8 to line -9 of page 2445
of the paper  that minimal radial functions of superharmonic
functions (with respect to a conformal Laplacian
of
a Riemannian metric)
are superharmonic.
Paper [TW] has made  essential use of the  false
assertion.
In our paper [LN]  we have used
at one point  $\min_{\partial B_r} v $, what 
Wang calls
a minimal radial function of $v$,
but we do not suppose that  $\min_{\partial B_r} v$
is superharmonic.

{\it So our objections have not been answered, and new
objections have arisen.
}

\section{}

In  the last paragraph of Section 3 of [W], Wang claimed
that Section 2 of [W] answers  our objections to the proof of the
main theorem  Theorem 1.3 in [TW].
  We disagree on that, and explain the reasons below.

Let us start by recalling Theorem 1.3 
in [TW]:

\bigskip

\noindent{\bf Theorem 1.3.} ([TW]) \ {\it  Assume that $\sigma$ satisfies 
${\mathbb C}_1-{\mathbb C}_4$, $\varphi\in C^0({\cal M})$,
$\varphi\ge c_0>0$, and $\Gamma$ is a convex cone
satisfying ${\mathbb  G}_1$ and ${\mathbb G}_2$.  Let $g_j=v_j^{  \frac 4{n-2} }
g_0$ be a sequence of solutions to (1.6).  Then
$v_j/\inf_{\cal M} v_j$ converges in $W^{1,p}$
(for any $1<p<\frac n{n-1}$) to an 
admissible function $v$.
Moreover, if $x_0$ is a singular point of $v$, then near $x_0$,
\[
v(x)=\frac {C_0+\circ(1) }{  d(x, x_0)^{n-2} }.
\tag{1.22}
\]
where $C_0$ is a positive constant, $d(x, x_0)$ denotes the geodesic 
distance from $x$ to $x_0$ in
the metric $g_0$.  Furthermore, each
singular point is isolated.
}

\bigskip

As defined in the first two lines of page 2443 of [TW], 
$x_0$ is a singular point of $v$ if there is a sequence of points
$x_j\in {\cal M}$ such that 
$v_j(x_j)/ \inf_{ \cal M}v_j\to \infty$ and $x_j\to x_0$.

In (1.22), the $\circ(1)$ term is to be understood in the 
usual sense, i.e. 
$$\lim_{ x\to x_0}
\left( d(x,x_0)^{n-2} v(x) -C_0\right)=0.
$$
This is agreed by Wang in the second to the last paragraph of 
Section 3 of [W].  Furthermore this usual sense of
convergence is    
 needed 
to  obtain Theorem 1.1 and 1.2 in [TW].
However we do not see this usual sense of
convergence in Theorem 1.3 is  being established in both [TW]
and Sections 2-3 of [W], as explained below.

The proof of Theorem 1.3 on page 2 of [W]
consists of three paragraphs.
We will call them paragraph 1,  paragraph 2, and  paragraph 3
respectively. 

We now phrase our question
as we follow these three paragraphs.
  We raise our question in the simplest situation
in order to more easily convey the ideas.

Paragraph 1 says: ``Let $x_j$ be the absolute maximum point of $v_j$.
As above we may assume 
 $x_j$ is a fixed point and $x_j=0$.
By Lemma 3.4,
the function $w_j$ given in
(3.11)   converges in
$W^{1,p}$ to the function $w$ in  
(3.29).
We need to show that $0$ is an isolated singular point of $w$.
''

Let us  look at
a simplified situation that 
$0$ is the only singular point of $w$.  More precisely,
$v_j/\inf_{\cal M} v_j$ is locally bounded in
${\cal M}\setminus \{0\}$ and
$v_j(0)/\inf_{\cal M} v_j\to \infty$.

 In this case, paragraph 2
is not needed, since that is used to prove that $0$ is an isolated 
singular point of $w$.
The part 
`` Therefore the absolute maximum point is an isolated singular point......
The above arguments also leads to a contradiction.'' 
in paragraph 3 is also not needed.
There is only one sentence left in paragraph 3 which
asserts that
the proof of Theorem 1.3 is completed.
{\it We do not see why this is the case --- estimate (1.22) 
has not been established}.

To make our point more precise, note that 
$$
w=- \frac 2{n-2} \log v,
\quad h(x):= w(x)-2\log |x|,
$$
so (1.22) is equivalent to
\begin{equation}
h(x)=\circ(1)\quad \mbox{in the usual sense as}\ x\to 0.
\label{usual2}
\end{equation}

\noindent However what 
 has been proved in Lemma 3.4 in the form stated in [W] is
(*) in [W], i.e.
$$
\lim_{ r\to 0}
\int_{ r<|x|<2r }
|h(x)|dx =0,
$$
which is much weaker than (\ref{usual2}).

\bigskip

\bigskip

[LN] Yanyan Li and Luc Nguyen,
{\it A compactness theorem for a fully nonlinear
Yamabe problem under a lower Ricci curvature bound},
arXiv:1212.0460v1 [math.AP] 3 Dec 2012.
\medskip

[TW] Neil Trudinger and Xu-Jia Wang,
{\it  The intermediate case
of the Yamabe problem for higher order curvatures},
International Mathematics Research Notices, Vol 2010,
no. 13, pp. 2437-2458.

\medskip

[W] Xu-Jia Wang,
{\it  Response to a question of Yanyan Li and Luc Nguyen
in their paper ``A compactness theorem for a fully
nonlinear Yamabe problem under a lower Ricci curvature bound'',
 arXiv:1212.0460},
arXiv:1212.3130v1 [math.AP].

\end{document}